\documentclass[twoside,leqno]{article}

\usepackage[letterpaper]{geometry}

\usepackage{siamproceedings}

\usepackage[T1]{fontenc}
\usepackage{amsfonts}
\usepackage{graphicx}
\usepackage{epstopdf}
\usepackage{enumitem}
\usepackage{algorithmic}
\ifpdf
  \DeclareGraphicsExtensions{.eps,.pdf,.png,.jpg}
\else
  \DeclareGraphicsExtensions{.eps}
\fi

\newsiamremark{remark}{Remark}
\newsiamremark{hypothesis}{Hypothesis}
\crefname{hypothesis}{Hypothesis}{Hypotheses}
\newsiamthm{claim}{Claim}
\newsiamthm{conjecture}{Conjecture}

\usepackage{amsopn}

\usepackage[utf8]{inputenc}

\usepackage{mathrsfs}

\usepackage{stix2}

\newcommand{\cc}{{\mathbb C}}
\newcommand{\pp}{{\mathbb P}}
\newcommand{\rr}{{\mathbb R}}
\newcommand{\qq}{{\mathbb Q}}
\newcommand{\zz}{{\mathbb Z}}
\newcommand{\nn}{{\mathbb N}}

\newcommand{\srD}{{\mathscr D}}

\begin{document}

\newcommand\relatedversion{}

\title{\Large Transcending Algebraic Geometry\relatedversion}
    \author{Carlos Simpson\thanks{C.N.R.S. UMR 7351, Laboratoire J.A.\ Dieudonné (\email{carlos.simpson@univ-cotedazur.fr}).}
    }

\date{}

\maketitle







\begin{abstract} 
We give an exposition of some of the ideas leading up to Jacob Tsimerman's groundbreaking work developing the theory of o-minimal definability in arithmetic and algebraic geometry. 
\end{abstract}

\section{Introduction}

Model theory has had a deep and extensive recent history of applications to geometry and arithmetic. The work of Jacob Tsimerman has brought together and pushed forward several strands of this line of thought, giving amazing new results such as: the André-Oort conjecture and the hyperbolic Ax-Lindemann theorem, first for the moduli ${\mathcal A}_g$ of abelian varieties and then for general Shimura varieties; their generalizations to the Ax-Schanuel theorem for $\zz$-variations of Hodge structure; definability of Griffiths period maps and the proof of the Griffiths conjecture that images of period maps are quasiprojective varieties; and the Shafarevich conjecture for quasiprojective varieties relative to linear representations of the fundamental group. 

This range of results gives us new insight on many facets of transcendental structures as they appear in algebraic geometry. The unifying tool is the use of the theory of {\em o-minimal structures}, and most particularly the o-minimal structure $\rr _{\rm an,exp}$ generated by analytic functions on bounded domains, and the real exponential function. 
Following precursors by Gabrièlov and Khovanskii, the study of these structures was initiated by van den Dries,  
followed by Pillay, Steinhorn and Knight, then Wilkie, and
continuing on the other hand with Macintyre, Marker, Miller, Denef, and others, as discussed for example in \cite{WilkieBourbaki,ScanlonBourbaki,PilaJT}. 

The theory of o-minimality as it can be applied to arithmetic geometry has been gaining steam ever since the introduction of model theory to geometrical questions by Hrushovski, and the consideration of generalized Zariski geometries by Hrushovski and Zilber. 
Some of the main questions to which this theory would apply came from works of Bombieri, Ax, André, Oort as shall be seen below. 

Since then, advances were made in stages by a number of mathematicians including Wilkie, Zannier, Pila,  Masser, Peterzil, Starchenko, Habegger, Klingler, Ullmo, Yafaev, Daw, Gao. Closely related are works of Yomdin, Edixhoven, Pink, to name a few. 

These advances, leading up to the many works that we shall be considering here by Tsimerman with Pila, Mok, Shankar, Bakker, Klingler, Brunebarbe, as well as to other generalizations that we'll mention briefly below, have notably been punctuated by a whole series of exposés in the {\em Séminaire Bourbaki} \cite{WilkieBourbaki,ScanlonBourbaki,FresanBourbaki}, and this theory has seen an exceptional number of publications in the {\em Annals of Mathematics}, {\em Inventiones} and other top journals. The mathematical community expressed a growing awareness of the transformative potential of this theory, but the consecration of these techniques as fundamental tools affecting all of geometry came only with Tsimerman's contributions. 

The flow starting from results on functional transcendence and leading towards entirely geometric applications in Hodge theory is a surprising and unexpected aspect of Tsimerman's work. 

In this note, we would like to emphasize how some of the most basic underpinnings of mathematics lead in a natural way to concepts such as transcendence and definability that are central to the theory.

\section{Beginnings}

As we learn Mathematics, its mysteries slowly enter our consciousness. It is worthwhile to look more closely at what happens at the very start of this process. 

Counting is not anything strange or wild, it is natural. On top of counting, we add layers of abstraction and new ideas that, when you think about it, are much more radical. 

The very first layer, a little outside of Mathematics and seen before we enter the Mathematical arena, is that we can write down symbols to express the results of counting. 
This idea of looking at the {\em expression} of a mathematical idea, is going to be fundamental to our exposition. 

The next layer is that these ``numbers'', whatever they are, can be manipulated by {\em operations}, and that these operation satisfy {\em rules} (commutativity, associativity, distributivity). And similarly, logical relations have operations ({\sc and}, {\sc or}, {\sc not}) that satisfy similar rules, and the rules themselves can be expressed with {\em quantifiers} $\forall , \exists$ 
that become, in turn, logical operations. 

Up until now, nothing too far out. The next stepping stone is the idea that there is 
an infinity of natural numbers $\nn$ (then, integers $\zz$, rational numbers $\qq$). This gives some consistence to $\forall$ and 
$\exists$ which start to say something that is not immediately perceptible in a finite way. And most of all, we can start to ask questions that do not necessarily have a cut-and-dried answer, an answer that can be attained by simple calculation using the operations. We can define prime numbers, make all sorts of computations, make conjectures, and prove some of them. 

The next rather strange construction is the set $\rr$ of real numbers. One already has several competing ways of defining these: decimal expansions (with a tricky equivalence relation saying that $1.0000\ldots  = 0.99999\ldots $), Cauchy sequences in $\qq$, Dedekind cuts. All of these concepts involve infinite collections of information, that is to say, functions defined on the infinity of natural numbers. 

\noindent
{\bf Who Shaves the Barber?} The set of real numbers is {\em uncountable}: there is no way of enumerating all the real numbers $a_1,a_2,a_3,\ldots $ in a list indexed by $\nn$. Indeed, if there were, we could define a new real number $\alpha$ whose $i$-th decimal, any one except $0$ or $9$ say, is different from the $i$-th decimal of $a_i$. Then $\alpha$ cannot be any of the $a_i$'s. 

The Constructivists say: wait a minute, let us just look at all the ways of writing down a formula for a real number. There are only countably many formulas, so there are really only countably many real numbers that we would ever need to worry about. This can be made into a valid theory (one never knows if the Barber is on the job today or not). One would call the resulting set of real numbers, the set of {\em constructible} or {\em definable} ones. Here we meet and use rather heavily concepts from the First Strange Beast, codifying what it means to have an expression that defines a real number. 

Let us work with the usual Cantorian idea of an uncountable set of real numbers where we do not require, at this stage, any constructibility properties. One can do a lot of new things with these real numbers, particularly, taking limits and hence infinite sums. In this way we can define many new real numbers, and whether or not one would like to look only at those that can be written down in some form, we certainly get {\em new} real numbers that were not rational, nor even solutions of polynomials with rational coefficients. 

Such things are called {\em transcendental}. They will be the subject of much of what we say below. 

On top of these ideas is our desire to draw pictures that represent things having to do with real numbers and their operations. For example, the graph of a function and other subsets of $\rr^2$. We then imagine being higher-dimensional beings that could draw pictures in $\rr^3$, $\rr^4$ and so on. Operations lead to subsets that can be pictured. 

As we progress further, we encounter more and more layers, such as the complex numbers which extend the real numbers in a convenient, although sometimes complicated way. But the ones listed above are the most basic ones leading to the rest. The amazing thing about mathematics is that we can actually say some things about how they behave. 

Jacob Tsimerman's works are among the most advanced and beautiful things we can say about the interactions between these basic ideas. 

\section{Algebraic definability}

If we restrict to the arithmetic operations: addition, multiplication, division by nonzero numbers, then we get the notion of algebraically definable subsets. Using the idea of a graph, we can look at definable subsets that could then be interpreted as graphs of functions. 

It turns out to be useful to use the {\em complex numbers} $\cc = \rr + \sqrt{-1} \rr$ here, since any nonconstant polynomial equation in complex numbers has a solution. 

The numbers definable using solutions of equations starting from the rationals form the field of algebraic numbers $\overline{\qq} \subset \cc$. These are the numbers that are contained in a finite degree field extension $K:\qq$. Since there are only countably many polynomials with rational coefficients, $\overline{\qq}$ is countable. The numbers that are not algebraic are called {\em transcendental}. 

As $\rr$ and hence $\cc$ are uncountable, there are infinitely many more 
transcendental complex numbers than algebraic ones. But, can we really write one down? Liouville gave a trick using decimal expansions with lots and lots of $0$'s. But we are interested in having expressions that allow to write down transcendental numbers, as we'll discuss shortly. But then, there are going to be only countably many such expressions, and if we reorder cleverly (part of the magic of $\nn$) one could say that ``every other such number is algebraic''. It is one of the paradoxes of such a discussion. 

But even if we manage to write down a new number, how can we know if it is transcendental? The study of transcendental numbers has a long and rich history. Some of the early results included proving that $\pi$ and $e$ are transcendental, then giving more general transcendence properties saying for example that numbers of the form $e^{\alpha _i}$ do not satisfy any non-obvious algebraic relations among themselves (when $\alpha_i$ are algebraic).

\section{Limits, sums and integrals}
\label{sec-limsum}

New kinds of real numbers may be obtained from the defining property of $\rr$, that it is the Cauchy completion of $\qq$. Thus, limits of sequances of rational numbers become real numbers that might well be transcendental. In practice, one meets such limits usually in the form of infinite sums, or integrals. For example the {\em exponential function}
is a sum
$$
e^x = \sum _{n=0}^{\infty} \frac{x^n}{n!} ,
$$
and its inverse the logarithm is an integral:
$$
\ln (y) = \int _{t=1}^{y} \frac{dt}{t}.
$$

Algebraic analysis has for two centuries or more studied functions defined by integrals starting from algebraically defined functions. These often satisfy differential equations, a phenomenon that is most often an expression of the Gauss-Manin connection that we'll meet later, in some way. 

{\bf Periods:} these are numbers obtained by integrating algebraically definable differential forms over domains of integration that represent homology classes in the algebraic varieties on which the integrands are defined and closed. 

The ``periods'' encompass many of the newly definable real or complex numbers that can be obtained by integration starting with algebraic functions. 

To continue a little with the history, early results included the {\bf Gelfond-Schneider theorem} saying that if 
$a,b\in \overline{\qq}$ are algebraic numbers, $a\neq 0,1$ and $b\not\in \qq$ then 
$a^b$ is transcendental. 
The {\bf Lindemann-Weierstrass theorem} gave algebraic independence of the numbers $e^{a_1},\ldots , e^{a_n}$ when $a_i\in \overline{\qq}$ are linearly independent over $\qq$.  In the 1960's, {\bf Baker's theorem} implied the generalization of Gelfond-Schneider to several numbers: if $a_i,b_i\in \overline{\qq}$, 
with $\log (a_i)$ independent over $\qq$ and  and $b_i$ not all rational, then  $a_1^{b_1} \cdots a_n^{b_n}$ is transcendental. 
Baker's theorem was stronger, giving a quantitative statement 
about relations among the logarithms.

These statements have formulations in terms of logarithms instead of exponentials, and one may view logarithms as basic kinds of periods. In this vein, other transcendence statements for periods, such as for the values of various functions associated to elliptic curves, enter into the classical repertoire \cite{Waldschmidt}. 

\section{Transcendence in algebraic and arithmetic geometry}
\label{sec-traag}

Periods have a rich history in algebraic geometry. They were main subjects of research by figures such as Euler, Gauss, Hermite, Abel, Jacobi, Riemann, Dedekind, Fricke, Klein, Picard, Fuchs, \ldots , leading into the development of wide swaths of algebraic geometry in the 20th century. 

Kontsevich and Zagier gave a general and simple definition of periods \cite{KontsevichZagier} as being any number defined as an integral of an 
algebraic function with $\qq$ coefficients, 
over a domain defined by polynomial inequalities with $\qq$ coefficients. These encompass what we usually think of as periods in algebraic geometry. Kontsevich and Zagier conjectured that any identities between periods should come from just a few simple operations: sums, change of variables, and Stokes formula. See the nice discussions in \cite{MullerStach,ViuSos} 
and the comparison with the Grothendieck period conjecture 
in \cite{Ayoub}.

The Kontsevich-Zagier period conjecture 
expresses what had been a general feeling for quite some time: that the new real numbers defined by integrating algebraic functions, should only obey identities that come about from a geometrical reason. 

One of the most famous incarnations of this idea is the {\em Hodge conjecture} that we'll explain next. Let us work with smooth projective algebraic varieties defined over subfields of $\cc$. A variety is a subset definable by polynomial equations in the complex projective space $\pp ^n_{\cc}$ that is the natural compactification of affine $n$-space $\cc^n$ obtained by adding the hyperplane $\pp ^{n-1}_{\cc}$ at infinity. Look at smooth varieties, those that are complex manifolds, these are the ones whose systems of defining equations satisfy the Jacobian condition to define smooth submanifolds at each point. 

The usual topology of $\cc^n = \rr ^{2n}$ extends to $\pp ^n_{\cc}$ and gives the {\em usual topology} of a smooth subvariety $X\subset \pp ^n_{\cc}$. This topological space, that we'll denote by $X^{\rm top}$, has its usual various topological invariants notably homology and cohomology. The algebraic variety structure leads to the {\em Hodge decomposition}
$$
H^i(X,\cc ) := H^i(X^{\rm top} ,\zz ) \otimes _{\zz } \cc 
= \bigoplus _{p+q=i} H^{p,q}(X).
$$
The sum takes place over positive values $p,q\geq 0$, so for example when $i=1$ we have just two terms: 
$$
H^1(X^{\rm top} , \cc )
= H^{1,0}(X) \oplus H^{0,1}(X).
$$
One furthermore has the condition $\overline{H^{p,q}} = H^{q,p}$. Indeed, the simplest way to define $H^{p,q}(X)$ is that it is the set of complex-valued cohomology classes that have representatives as de Rham cohomology classes which are harmonic forms containing $p$ times coordinate forms $dz_j$ and $q$ times coordinate forms $d\overline{z}_k$ at
each point. This characterization, and the direct sum property, use the property that when $X$ is viewed as a complex manifold, it has the additional strong property of being a {\em Kähler manifold}. 

The {\em rational Hodge classes} (in even dimensional cohomology) 
are the elements of the intersection, taken within 
$H^{2p}(X^{\rm top} ,\cc )$:
$$
{\rm Hdg}^{2p}(X,\qq ):= 
H^{2p}(X^{\rm top} ,\qq ) 
\cap H^{p,p}(X).
$$
A basic observation is that if $Y\subset X$ is a codimension $p$ subvariety (it doesn't need to be smooth) then its fundamental class $[Y]$ is a Hodge class (the degree of this class is the real codimension of $Y$, namely twice the complex codimension $p$). 

The {\bf Hodge conjecture} states the converse: {\em all the Hodge classes come from fundamental classes of subvarieties}. An important variation on this statement by André replaces algebraic cycles with {\em motivated cycles} obtained by formally adjoining the adjoint of the Kähler operator to the algebraic cycles \cite{Andre}. 

The {\em Hodge filtration} is defined as the sum of the Hodge subspaces going in one direction, $F^pH^i := \bigoplus _{r\geq p} H^{r,i-r}(X)$. These subspaces 
\begin{itemize}

\item define the Hodge decomposition because of the complex conjugation relation; 

\item and vary in a holomorphic way when the equations of $X$ are moved, indeed they are defined by period integrals. 

\end{itemize}

It is instructive, and will be useful for our main discussion later, to consider the case of $H^1$. In this case, the space of regular algebraically defined holomorphic $1$-forms on $X$, denoted $H^0(X,\Omega ^1_X)$, defines by integration a collection of periods:
$$
H^0(X,\Omega ^1_X) \times H_1(X^{\rm top},\zz ) \rightarrow \cc , 
\;\;\;\; (\alpha , \xi ) \mapsto \int _\xi \alpha  
$$
hence giving a subspace of cohomology 
$$
H^0(X,\Omega ^1_X) \hookrightarrow H^1(X^{\rm top} ,\cc ).
$$
This subspace $H^{1,0}(X)$. 

More generally, all the subspaces of the Hodge filtration may be defined as periods by looking at certain meromorphic differentials that represent hypercohomology classes of the truncated de Rham complex. 

The {\em Hodge structure} of the variety is defined as the triple 
$(V,V_{\zz}, \{ V^{p,q} \} )$ of a complex vector space, an integer lattice, 
and a decomposition satisfying $\overline{V^{p,q}} = V^{q,p}$. it encodes 
the periods of certain algebraic differentials on $X$, in a way most closely related to the topology and geometry of $X$. 

At its most basic level, a Hodge structure may be pictured as an integer lattice together with various subspaces, and the coordinates of these subspaces with respect to the lattice are periods---in general, transcendental numbers. 

\section{Siegel's upper half space and Griffiths period domains}

Griffiths made an extensive and deep study of the spaces, called {\em period domains}
\cite{CarlsonGriffiths}, that parametrize Hodge structures, and how the Hodge structures of a family of varieties define the corresponding map from the base space of a family into the period domain. 

A good collection of articles, containing both introductory expositions and advanced properties, is \cite{transag}. 

Historically, the first example in higher dimensions was the {\em Siegel upper half space} parametrizing abelian varieties. We make a small digression to talk about  abelian varieties. These are smooth projective varieties that have an abelian group law
$$
A\times A \rightarrow A.
$$
There is what we shall be calling the {\em Lie group viewpoint}, that the universal covering is a complex vector space and the fundamental group is an integer lattice:
$$
\widetilde{A} \cong \cc ^g, \;\;\;\; \pi _1(A) \cong \zz ^{2g}, \;\;\;\; , 
A^{\rm an} \cong \cc ^g / \zz ^{2g}.
$$
In the presence of a {\em polarization form} (antisymmetric nondegenerate bilinear form on $\zz^{2g}$ such that $\cc^g$ becomes an isotropic or Lagrangian subspace with positive definite metric) such a 
``complex torus'' is in fact a smooth projective variety. The varieties that arise this way are recognized as those admitting a group structure (necessarily abelian). It follows that being an abelian variety is a property that can be defined in terms of the equations defining $A$. The viewpoint of polynomial defining equations, in principle for an embedding $A\hookrightarrow \pp^N_{\cc}$, will be called the {\em algebraic geometry viewpoint}. 

Going back to the Lie group viewpoint, it is not hard to see that the Hodge structure on $H^1(A^{\rm top}, \cc ) \cong \cc ^{2g}$ determines the {\em period matrix} 
$$
\Phi : \zz ^{2g} \rightarrow \cc ^g
$$
that in turn determines the structure of $A$ as a complex torus. This in turn leads to the {\em theta functions} (some complicated sums of functions over $\cc^g$) that give an embedding $A\hookrightarrow \pp ^N_{\cc}$ getting us back to the algebraic geometry viewpoint. 

We now see a path to parametrizing abelian varieties, in the Lie algebra viewpoint: if we can look at the appropriate subspace of the Grassmanian $G(\cc^{2g},\cc ^g)$ corresponding to the polarization conditions, then the quotient of this subspace by the group representing change-of-basis of $\zz^g$ will parametrize the abelian varieties. 

The previous paragraph exactly describes the {\em Siegel upper half-space} 
${\mathbb H}_g$ with an action of the group $Sp(2g, \zz )$ such that the quotient is the moduli space ${\mathcal A}_g$
of abelian varieties. This is the moduli space, viewed itself in the Lie group perspective. 
It may be written as 
$$
{\mathbb H}_g = Sp(2g,\rr )  / U(g), \;\;\;\;
{\mathcal A}_g = Sp(2g,\zz )  \backslash {\mathbb H}_g  = 
 Sp(2g,\zz )\backslash Sp(2g,\rr ) / U(g) .
$$
The concept of a {\em period matrix} also allows us to write an expression for 
${\mathbb H}_g$ that looks very much like a matrix version of the usual upper half plane. The point in the Grassmanian $G(\cc^{2g},\cc ^g)$ may be normalized, by changing basis in $\cc^g$, to come from a $2g\times g$ matrix of the form 
$$
P = \left[ \, Id_{g\times g} \,  | \, Z \, \right] 
$$
where $Z \in Sym_{g\times g} (\cc )$ is a $g\times g$ symmetric matrix. 
The {\em polarization condition}
is  $Im(Z ) > 0$ (meaning that this real symmetric matrix is positive definite). 
The action of the symplectic group $Sp(2g,\zz )$ may be written in the form 
$$
\left( \begin{array}{cc}
A & B \\
C & D 
\end{array} \right) \cdot Z := 
(A Z + B) (CZ + D) ^{-1}
$$
which looks a lot like the matrix version of the usual action of $SL(2,\zz )$ on the 
classical upper half plane. 

One of the first astounding results in this topic is the theorem that the above Lie-algebra picture has a corresponding algebraic geometry picture.

\begin{theorem}[Satake,Baily-Borel]
This quotient admits a compactification $\overline{{\mathbb H}}_g$ that is a projective algebraic variety\footnote{More precisely, it is a Deligne-Mumford stack; we shall be systematically ignoring the stacky aspects of these moduli spaces.} containing 
${\mathbb H}_g$ as a Zariski-open subset. Furthermore, these varieties are defined over $\qq$. 
\end{theorem}

Mumford's GIT provided a very general framework in which to make such moduli constructions, and allowed to provide a construction of the moduli space, valid over arbitrary ground rings, not depending on the Lie group viewpoint. 

More general {\em Shimura varieties} are obtained by replacing ${\mathbb H}_g$
with a general Hermitian symmetric space, and their compactifications and algebraic geometry embeddings are also known. 

We meet here an important difference with the more general Griffiths period domains. Griffiths defined a {\em transversality condition} for maps into a period domain. When it is Hermitian symmetric, the transversality condition is vacuous, and all points have a meaning in algebraic geometry. However, in the general case, basically for higher weight Hodge structures, the transversality condition is a nontrivial constraint on any map into the period domain that comes from algebraic geometry. It follows that most points of the period domain will not be related to algebraic geometry, and the development of an adequate ``algebraic geometry'' viewpoint is much more delicate. That is indeed the topic of some of Tsimerman's most recent contributions. 

For most of the rest of the present exposition, we are going to look at the case of the moduli of abelian varieties which is the theater of one of Tsimerman's early signature works. 

\section{Special subvarieties I: the Manin-Mumford and Mordell-Lang conjectures}

Next we are going to look at a ``baby case'' of the theorem about special points and special subvarieties. Consider an abelian variety $A$, presented in terms of equations as $A\subset \pp^N_{\cc}$ (let's say for example that the equations have coefficients in a subfield $K\subset \cc$) , and presented in the Lie-theory viewpoint as 
$A=\cc^g / \Lambda$ where $\Lambda \cong \zz ^{2g}$ is an integer lattice. 

One may define some {\em special points}, namely the {\em torsion points}. The group law $m: A\times A \rightarrow A$ is an algebraically defined map in the algebraic geometry viewpoint, but it is also just addition of vectors in $\cc^g$ in the Lie theory viewpoint. Let $o\in A$ denote the origin corresponding to the equivalence class of $0\in \cc^g$. A point $x$ is {\em $m$-torsion} if $x+x+\cdots x = o$ (the sum of $m$ elements). This condition is defined by algebraic equations in our algebraic viewpoint, so for example the torsion points will have coordinates in $\overline{K}$. But it is also the condition that $x$ is the image of a point $y\in \cc^g$ such that $m\cdot y \in \Lambda$. In other words, the set of  $m$-torsion points take the form 
$$
A_{{\rm tors},m} = \left( {\textstyle \frac{1}{m}}\Lambda \right) / \Lambda \subset \cc ^g / \Lambda .
$$
The coordinates of $y\in \cc^g$ will therefore be rational multiples of the periods i.e. the coordinates in $\cc^g$ of the basis vectors of $\Lambda$. Denote by $A_{\rm tors}$ the union of the subsets of $m$-torsion points for all $m$. 

We also have some {\em special subvarieties}: these are the sub-abelian varieties $B\subset A$, that is to say the subvarieties defined by algebraic equations that are compatible with the group law. In the Lie theory viewpoint, these correspond to sub-Hodge structures, that is to say subspaces $\cc ^d \subset \cc^g$ such that $\Lambda ' = \cc ^d \cap \Lambda$ is a lattice of rank $2d$. In either viewpoint, this is not an easy condition to satisfy, and for ``general points'' (think of ``choosing an abelian variety at random'') there will not be any sub-abelian varieties. However, they can exist for example in the case of products 

We may combine these two notions of speciality and look at special subvarieties that are translates of sub-abelian varieties by torsion points. 

\begin{theorem}[Manin-Mumford conjecture] 
\label{mmconj}
Suppose $C$ is an algebraic curve of genus $\geq 2$, and let $J$ denote its {\em Jacobian abelian variety} with embedding $\alpha : C \hookrightarrow J$. Then $\alpha (C) \cap J_{\rm tors}$ is finite. 
\end{theorem}

This was conjectured by Manin and Mumford and proven by Raynaud. Serge Lang formulated a more general statement viewed as a generalization of the Mordell conjecture (Faltings' theorem), and his formulation meets most closely with our comparison between the Lie theory and algebraic geometry viewpoints. 

\begin{theorem}[Mordell-Lang conjecture]
\label{mlconj}
Suppose $Z\subset A$ is a subvariety in the algebraic viewpoint, that is a closed subset defined by complex algebraic equations. 
Suppose $\Lambda \subset A$ is a subgroup of finite $\qq$-rank, and suppose $Z\cap \Lambda$ is {\em Zariski-dense} in $Z$,
that is to say not contained itself in a smaller closed algebraic subvariety. Then $Z$ is a translate of an abelian subvariety. 
\end{theorem}

In the case $\Lambda = A_{\rm tors}$ it was also first proven by Raynaud \cite{Raynaud}, following a similar result by Bogomolov \cite{Bogomolov}, and 
the proof for   
general $\Lambda$ in characteristic $0$ was due to Faltings \cite{Faltings}, 
with Hindry  removing the hypothesis that $\Lambda$ be finitely generated \cite{Hindry}.
It was reproved and
generalized by many people including Vojta and McQuillan for semi-abelian varieties 
\cite{McQuillan}. A new proof of \ref{mmconj} 
and generalization of \ref{mlconj} 
 to the finite characteristic function field case by Hrushovski \cite{HrushovskiManMum,Hrushovski} was one of the original applications of model theory to arithmetic geometry. 

Let us contemplate what this basically says: that any polynomial identity that holds, in the algebraic geometry picture, for many special points at a time, has to be of a known form coming from a (quasi-)special
subvariety where we can see easily that the identity should hold in terms of the Lie theory picture. 

The full statement fits into a somewhat more general formulation that can be summed up as follows.

\begin{definition}
\label{mmmlao}
Suppose given a variety $M$ with a notion of {\em special subvariety} and a subset $\Lambda \subset M$ of {\em special points}. We say that the {\em MMMLAO property} (for Manin-Mumford-Mordell-Lang-André-Oort) holds if there is a finite union of special subvarieties $M^{\rm spe} \subset M$ such that $M^{\rm spe}\cap \Lambda$ is Zariski-dense in $M^{\rm spe}$, and $(M-M^{\rm spe})\cap \Lambda$ is finite. 
\end{definition}

The Mordell-Lang conjecture gives this property for the abelian variety with special subvarieties being translates of sub-abelian varieties and special points being the torsion 
points. We'll see, along the exposition, how this property interplays with other related properties. 

\section{Special subvarieties II: CM points and Mumford-Tate groups}

In our prior discussion about the moduli space of abelian varieties, we have bumped up a level to go from looking inside a given abelian variety, to looking at the moduli space parametrizing the collection of all abelian varieties as a whole. In the bumped-up version, we also have Lie theory and algebraic geometry viewpoints, with the added bonus that the Lie groups showing up can now be non-abelian such as the symplectic group $Sp(2g)$. 

It is natural to try to extend the concepts of special points and special subvarieties to this case. How are we going to do this? The idea will be to relate properties happening in the Lie theory viewpoint, with properties happening in the algebraic geometry viewpoint. We have a virtual tool to do that: the Hodge conjecture. 

Even though its status remains highly unknown and open, the Hodge conjecture leads to specific predictions about the behavior of special relationships between periods, that have been proven. Suppose we are given a family of varieties $f:X \rightarrow S$. It means an algebraic morphism whose fibers are smooth and projective, that we roughly think of as creating a family of varieties $X_s:= f^{-1}(s)$ that vary as a function of the point $s\in S$. 

A family of varieties then leads to a family of Hodge structures
$$
S\ni s \mapsto \left( H^i(X_s^{\rm top},\cc ), H^i(X_s^{\rm top},\zz ),
\{ H^{p,q}(X_s) \} \right) .
$$
This family has a {\em Gauss-Manin connection} reflecting the transportation of cohomology classes from $X_s^{\rm top}$ to $X_t^{\rm top}$ along any path
from $s$ to $t$, but the subspaces $H^{p,q}(X_s)$ vary in a way that satisfies Griffiths' axioms for a {\em variation of Hodge structure}. 

In this situation, the Hodge conjecture leads us to look at the set of points $s\in S$ where $X_s$ acquires an additional Hodge cycle $\eta \in H^{2p}(X_s^{\rm top},\qq )
\cap H^{p,p}(X_s)$. One can define precisely the irreducible components of the space parametrizing ``extra Hodge cycles'', known as the {\em Noether-Lefschetz locus}. 

The Hodge conjecture predicts that Hodge cycles come from algebraic cycles in the fiber algebraic varieties. The condition of having an algebraic cycle is an algebraic-geometry one, that can be expressed in terms of equations. (Some additional acrobatics here using the algebraic de Rham cohomology expression for $H^i(X_s^{\rm top},\cc )$ is useful here.) If the Hodge conjecture were known to be true, then, the Noether-Lefschetz locus would be definable in terms of algebraic equations. 

This turns out to be a theorem proven by Cattani, Deligne and Kaplan  \cite{CattaniDeligneKaplan} based on the huge body of work by Griffiths, Schmid, Cattani, Kaplan, Kashiwara, Kawai and others on the local properties of variations of Hodge structure. Thus, even though the Hodge conjecture is not known, its consequence saying that the Noether-Lefschetz locus is defined algebraically, is a theorem. An interesting side note is that Bakker, Klingler and Tsimerman have given a new proof of the result of \cite{CattaniDeligneKaplan}
using o-minimal structures \cite{BKT}. 

Of course, having a Hodge cycle is a relatively simple linear algebra condition in the Lie theory viewpoint on Hodge structures. It therefore satisfies our criterion for what we would like to use to define a notion of {\em special subvariety}. In a nutshell, given a family of varieties $f:X\rightarrow S$, the components of the Noether-Lefschetz locus will be considered as the special subvarieties of $S$.

This applies to our moduli space ${\mathcal A}_g$ of abelian varieties. We have a universal family\footnote{Ignoring, as we are doing, the stacky aspects of this question: to be correct one should use the Deligne-Mumford stack version of ${\mathcal A}_g$ as the base of the universal family.}
$f:{\mathcal X}_g \rightarrow {\mathcal A}_g$, and we may say that the {\em special subvarieties} of $ {\mathcal A}_g$ are those loci where ``extra Hodge cycles'' are acquired. Let us note another {\em caveat}: one needs to consider Hodge cycles in the tensor algebra of $H^1(A_s)$, indeed Hodge cycles only occur in even dimensions so $H^1(A_s)$ itself doesn't have any. The theory adapts easily to this extension. 

The {\em Mumford-Tate group} provides a convenient way of organizing the information about additional Hodge cycles and Noether-Lefschetz loci. One way of formulating it is by the theory of Tannakian duality due to Grothendieck and Saavedra Rivano, see Deligne's exposition\footnote{In \cite{LNM900},
Deligne also proves the property ``Hodge implies absolute Hodge'' for abelian varieties, one of the first that went in the direction of showing a consequence of the Hodge conjecture.} 
in \cite{LNM900}. Given a Hodge structure, such as $H^1(A )$ for an abelian variety $A$, we can form the Tannakian category of all subquotient Hodge structures in the tensor algebra (including tensors with the dual). The {\em Mumford-Tate group} is the algebraic group corresponding to this Tannakian category by Tannaka duality. 

It may also be characterized as the subgroup $G_{MT} \subset GL(H^1(A,\qq ))$
of elements that preserve all of the Hodge cycles in the tensor algebra. 

The {\em special subvarieties} of ${\mathcal A}_g$ are now defined as the subvarieties of points where the Mumford-Tate group jumps down to a new lower size. 

Defined in this way, it now contains arithmetic information too: there are components of the special locus of dimension $0$, i.e. special points. These are the points where $G_{MT}$ is the smallest possible. For the case of abelian varieties, the structure of the resulting points can be analyzed precisely. The special points are called {\em CM-points}, with ``CM'' standing for Complex Multiplication, by analogy with the original case of elliptic curves. 

An abelian variety is of CM type if its Hodge structure splits, when the rational vector space is tensored by an appropriate field extension, into pieces of rank $1$. This means that there is a field $K$ such that $H^1(A,\qq )\otimes _{\qq} K = \bigoplus V_i$
with $V_i$ being rank $1$ $K$-vector spaces, such that when we further tensor from $K$ to $C$, the $V_i$ group together in two parts to form the $H^{1,0}$ and $H^{0,1}$ of the Hodge decomposition. In sum, it is the easiest way that we would be able to write down explicitly the Hodge decomposition. 

A little bit of linear algebra shows that the field $K$ turns out, in this case, to also be the endomorphism algebra (tensored with $\qq$) of the Hodge structure $H^1(A)$. And, since the Hodge structure determines the complex torus which in turn determines the subvariety $A\subset \pp ^N_{\cc}$ in the algebraic geometry viewpoint, we find that a subring
generating the field $K$ appears as the algebra of endomorphisms of the algebraic-variety-with-group-structure $(A,+)$. 

The CM condition therefore also has an expression in terms of algebraic equations defining points of the moduli space ${\mathcal A}_g$. We call the special points defined in this way the {\bf CM-points} of  ${\mathcal A}_g$.

\section{Results}

In this section we consider the moduli space of abelian varieties ${\mathcal A}_g$ which is, in some sense, the main irreducible Shimura variety in higher dimensions. See Section \ref{sec-shimura} for statements in the general setting. 

In ${\mathcal A}_g$, the set of CM points is completely analogous to the set of torsion points in a single abelian variety, and it is natural to as whether the Manin-Mumford-Mordell-Lang theorem has a generalization to this setting. 

This is Tsimerman's theorem \cite{TsimermanAnnals}:  

\begin{theorem}[Tsimerman: the André-Oort conjecture for moduli of abelian varieties] 
\label{andreoort}
Suppose $Z\subset {\mathcal A}_g$
is an algebraically defined subvariety such that the set of CM-points contained in $Z$, is Zariski-dense in $Z$. Then $Z$ is a special subvariety defined by jumping of the Mumford-Tate group. 
\end{theorem}

Let us separate out one of the theorems that plays a central role in the proof. It is itself a natural statement that the reader may well have already wondered about. We recall that Griffiths period domains are usual open subsets of their corresponding flag varieties, since there is a positive definiteness condition on the bilinear form induced by the  polarization for the Hodge subspaces. A close cousin is the unit ball contained in projective space, that is indeed a Shimura variety for the group $PU(n,1)$. 

Thus in order to talk about algebraic definitions of subsets of period domains (including of course our main ${\mathbb H}_g$) we should speak of {\em semi-algebraic sets}, these are subsets defined by polynomial equations and inequalities. 
It is now natural to ask, how transcendental is the projection 
${\mathbb H}_g\rightarrow {\mathcal A}_g$ in terms of algebraically defined subsets on both sides? This is the topic of the {\em hyperbolic Ax-Lindemann theorem} (based on Ax's concept of functional transcendence \cite{Ax}) proven by Pila and Tsimerman \cite{PilaTsimerman}: 

\begin{theorem}[Pila-Tsimerman: hyperbolic Ax-Lindemann for moduli of abelian varieties]
\label{axlindemann}
Suppose $V\subset {\mathcal A}_g$ is an algebraic subvariety whose inverse image in 
${\mathbb H}_g$ contains a 
semi-algebraic set $W$. Then there is a weakly special subvariety $S\subset V$ containing the image of $W$. In particular \cite{UllmoYafaev}, if the full inverse image is algebraic in ${\mathbb H}_g$ then $V$ was weakly special. 
\end{theorem}

The notion of ``weakly special'' here is the appropriate analogue in the present situation of the notion of an arbitrary, rather than torsion, translate of a sub-abelian variety in the setup of Mordell-Lang, see \cite[\S 2.4]{PilaTsimerman}. 

One may in fact formulate an Ax-Lindemann theorem inside abelian varieties themselves, with a statement that the reader may intuit from the above: if we have an algebraic subvariety of $A$ whose inverse image by the projection 
$\cc^g\rightarrow A$ is also an algebraic subset in $\cc^g$, then the subvariety has to be a translate of a sub-abelian variety. This entered into the work of Pila-Zannier giving a new proof of Manin-Mumford-Mordell-Lang \cite{PilaZannier,Pila}. See the exposition by Orr \cite{Orr} for a discussion as well as a useful reference, see also Scanlon's Bourbaki \cite{ScanlonBourbaki}.

\section{Inputs from number theory}

Most of the article \cite{TsimermanAnnals} is devoted to the number-theory results and arguments needed to approach the situation of the Ax-Lindemann theorem \cite{PilaTsimerman,KlinglerUllmoYafaev} that is then applied to obtain special subvarieties. 

We will try in this section to give an idea of how and why these arguments work. 
The techniques and steps are sometimes referred to as the {\em Zannier strategy}, introduced in the article of Pila and Zannier \cite{PilaZannier}, and subsequently applied by Pila \cite{Pila} as well as in numerous other places. We refer to \cite{Pila,ScanlonBourbaki,DawSurvey,PilaJT} for further references.

\subsection{Heights}
In number theory, the {\em height} of a number in $\overline{\qq}$ is a quantity measuring mainly the size of its denominator, but also the size of the field extension $K:\qq$ in which it lies. This extends in a natural way to a measure of points in an algebraic variety defined over some algebraic number field. 

In the case of a point in ${\mathcal A}_g$ that parametrizes an abelian variety $A$, Faltings defined a normalized height that is called the {\em Falting height} $h_{\rm Falt}(A)$ with good invariance properties. 

\subsection{The Colmez conjecture}
The points we are looking at are CM abelian varieties, whose isogeny classes are determined by the CM field $E$ (it should satisfy appropriate conditions: being a totally imaginary extension of a totally real field) and a choice $\Phi$ of ``CM-type'' which means choosing half of the embeddings $E\hookrightarrow \cc$, the other half given by complex conjugates. 

One main ingredient is a statement conjectured by Colmez, stating that the Faltings height of a CM abelian variety depends only on $(E,\Phi )$ and giving a 
precise formula with good asymptotics. 

While this remains open, Tsimerman manages to use instead an averaged version of the Colmez conjecture that has been proven by 
Andreatta, Goren, Howard and Madapusi-Pera \cite{AGHMP} and 
Yuan and Zhang \cite{YZ}. 

The consequence of this result is that the Faltings height is bounded sub-polynomially by the discriminant of the CM field $E$ \cite[Corollary 3.3]{TsimermanAnnals}. 

\subsection{The Pila-Wilkie counting theory}

The concepts of model theory, definability and $o$-minimal structures enter into the proof of Theorem \ref{andreoort} in essentially two different places. We are going to give a general discussion of the theory below, concentrating on how it applies to the hyperbolic Ax-Lindemann theorem as well as related and further directions. 

For now, we need the notion of definability in the number theoretical part in order to state the Pila-Wilkie theorem \cite{PilaWilkie}. See Section \ref{sec-ranexp} below for the definition of
$\rr _{\rm an,exp}$. The map $\pi : {\mathbb H}_g \rightarrow {\mathcal A}_g$ is definable, when restricted to the closure of a fundamental domain.

\begin{theorem}[Pila-Wilkie]
\label{counting}
Suppose $Y\subset \rr^n$ is a definable subset in an $o$-minimal structure such as $\rr _{\rm an,exp}$. If $Y\cap \qq ^n$ has a lot of points in terms of the height, meaning that the number of points of height $\leq H$ grows as a strictly positive power of $H$, then $Y$ has to contain a positive-dimensional semialgebraic subset. 
\end{theorem} 

See \cite{PilaWilkie} for a discussion of the history of this result, going back for example to Sarnak \cite{Sarnak} and Yomdin \cite{Yomdin}. It is a key step in the approaches of many people to proving variants of the Manin-Mumford and Mordell-Lang conjectures. 
In further work, a generalization to the non-archimedean case was done by Cluckers, Comte and Loeser \cite{CCL}, and a conjecture of Wilkie was proven by Binyamini-Novikov \cite{BinyaminiNovikov} and Binyamini-Novikov-Zak \cite{BinyaminiNovikovZak}.

A couple of important points may be illustrated by some examples. First, one might ask whether to expect a stronger MMMLAO type of property (Definition \ref{mmmlao}) to hold in this setting. Namely, could there be a result saying that for a definable subset $Y$, if we define $Y^{\rm alg}$ to be the union of the positive dimensional algebraic subsets of $Y$, then will $(Y-Y^{\rm alg})\cap \qq^n$ be finite? The following example shows that this is not the case: let 
$$
\rr^2 \supset Y := \{ (x, 2^x) \} 
$$
be the graph of the function $x\mapsto {\rm exp}((\ln 2)x)$. This is clearly $\rr _{\rm an,exp}$-definable, but it contains the infinite collection of points of the form $(k,2^k)\in \zz^2$ for $k\in \zz$. We notice here that the heights of these points grow exponentially, so the number of points of height $\leq H$ grows as $\log H$ as stated in the counting theorem. The Pila-Wilkie counting property is the appropriate statement in this setting, analogous to the MMMLAO property but with ``finite'' replaced by the sub-polynomial bound for the number of points of a given height. 

The second example, very similar to the first, is to let $Y$ be the graph of the function $x\mapsto \sin (\pi x)$. This contains the points $(k,0)\in \zz^2$. In this case, the height of the $k$-th point is just $k$ so the number of points of height $\leq H$ grows as $H$, a strictly positive power. This contradicts the bound of the counting theorem, but in this case the subset $Y$ is not definable in any o-minimal structure. Indeed, the intersection of $Y$ with the $x$-axis has infinitely many distinct points, which is not allowed in the definition of an o-minimal structure. This shows that the o-minimal definability is a crucial part of the counting theorem. 

\subsection{The proof of Theorem \ref{andreoort}}

The Pila-Wilkie theorem is the step leading to a semialgebraic subset of ${\mathbb H}_g$. Tsimerman does this by applying the averaged Colmez conjecture mentioned in the previous subsection, to get enough CM points whose height, in the Lie theory picture in 
${\mathbb H}_g$, is bounded. 

This is based on a beautiful argument using Galois orbits and the field of definition of the subvariety of ${\mathcal A}_g$ that we are looking at. This part is the main new ingredient allowing Tsimerman to go to the proof of the André-Oort property, but it would go beyond our present scope (and the competency of the author) to discuss the details. 

The Ax-Lindemann theorem \ref{axlindemann} permits to go from a semialgebraic subset of ${\mathbb H}_g$ to a weakly special subvariety, which is indeed special since it contains special points. 

\section{Model theory and definability structures}

Here starts the second part of our exposition, devoted to the theory of definable sets in an o-minimal structure and leading towards the applications in Hodge theory. The reader has been encouraged all along to think about the idea of how we define subsets in terms of algebraic equations and other constructions,
functions and inequalities. These natural mathematical thoughts are made precise in mathematical logic in the subject of {\em model theory} with its core notion of {\em definability}. 

One of the first to apply model theory to questions in number theory and geometry was Hrushovski. This led to a whole body of work by many mathematicians, leading up to the many applications of the notion of definability that we are considering here. 

A {\em predicate} $P(x_1,\ldots , x_n)$ is a function, depending on variables $x_1,\ldots , x_n$, taking truth values, that is to say taking values in the most basic boolean algebra 
$\{ {\rm True}, {\rm False} \}$. One can make important distinctions between various versions of this notion in type theory. For us, an important idea is that a predicate can sometimes be written down in more or less explicit terms. 

Giving a predicate is equivalent to giving the subset $\Sigma _P$ of values of the 
$n$-tuple $(x_1,\ldots , x_n)$ on which $P$ is ${\rm True}$. For us, the variables are going to take real number values, so the full set of $n$-tuples $(x_1,\ldots , x_n)$
is $\rr^n$ and giving the predicate $P(x_1,\ldots , x_n)$ is equivalent to giving the subset $\Sigma _P \subset \rr ^n$. 

Boolean operations on predicates correspond to set theoretical operations on subsets: 
$\cap$ for ${\rm and}$, $\cup$ for ${\rm or}$, and complementation for negation. Furthermore, quantifiers take on an easy expression in these terms: if we let 
$$
Q(x_1,\ldots , x_{n-1}) := \exists x_n, P(x_1,\ldots ,x_n)
$$
then $\Sigma _Q \subset \rr ^{n-1}$ is just the image of $\Sigma _P \subset \rr^n$
by the last projection $pr_n:\rr ^n \rightarrow \rr^{n-1}$, and similarly of course for any of the other coordinates. The quantifier $\forall$ is obtained by complementation. 

The definability structure $\rr _{\rm alg}$ consists of all the subsets, or predicates, that can be defined using arithmetical operations, the order relation $<$, and these boolean quantifiers on the variables. This is the most basic notion of definability, that we have indeed been using all along in the concept of algebraic 
equations. 

This mirrors closely Chevalley's notion of {\em constructible set} in algebraic geometry, with the difference that the order relation $<$ is now included whereas for constructible sets only equalities and inequalities of the defining polynomial relations were used. 

More generally a {\em definability structure} is a collection of subsets $\Sigma \subset \rr^n$ that are called the {\em definable sets} for that structure, that is closed under boolean operations, pullbacks, and quantifiers (in other words, under images of projections). The {\em definable functions} are then the functions such that the graph is a definable set (in particular the domain is a definable set). 

\subsection{o-minimality}
There is a remarkable property that definability structures can have, called {\em o-minimality}. A structure is said to be {\bf o-minimal} if the definable subsets of $\rr$ are finite unions of ``intervals'', where an {\em interval} is defined as $\rr$ or an intersection of open or closed half lines i.e. $(a,\infty )$, $[a,\infty )$, $(-\infty , a)$, or $(-\infty , a]$. 
Notice that $\emptyset$, singleton points, and intervals from $a$ to $b$ possibly including endpoints, are included. 

This notion is a fulcrum between examples and properties. It is amazing that there exist nontrivial examples of o-minimal structures, including mainly $\rr _{\rm an,exp}$ as we shall see below, but also that definabiity in an o-minimal structure leads to many highly advanced properties that come into play in the theories we are discussing. 

Let us notice some functions that are {\em not} definable in an o-minimal structure. For example, $\sin (x)$, $\cos (x)$ and more generally pretty much any periodic function is not going to be definable, otherwise the intersection of the graph with a  horizontal line would define a subset of $\rr$ that is a countably infinite collection of points, contradicting the o-minimality property.

\subsection{The structure $\rr _{\rm an,exp}$}
\label{sec-ranexp}
Let us look now at the definition of the main o-minimal structure used in Tsimerman's theories. We define the structure $\rr _{\rm an,exp}$ as an extension of $\rr _{\rm alg}$ defined by adding certain functions. The first are the {\bf real analytic functions} on bounded domains. It is sufficient to add the real analytic functions $f(x_1,\ldots , x_n)$
with domain of definition a box $|x_i - a_i| \leq c$ such that $f$ admits an extension to a real analytic function on a slightly bigger box $|x_i - a_i| < c+\epsilon$.

Real analytic functions on more general relatively compact domains are obtained by taking the unions of graphs, in order to glue together functions, then intersecting with the pullback of the desired domain (that should itself be definable). 

After the real analytic functions on relatively bounded domains, we add the real {\bf exponential function} $x\mapsto e^x$. Then $\rr _{\rm an,exp}$ is completed to the smallest definability structure containing these functions. 

The o-minimality of the $\rr _{\rm an,exp}$ structure is due to 
van den Dries and Miller \cite{vdDriesMiller} with subsequent investigations by
van den Dries-Macintyre-Marker \cite{vdDriesMacintyreMarker} and then leading to a wide study by many people, as may be seen in the numerous references. 
We note in particular the series of articles by Pillay, Steinhorn \cite{PillaySteinhorn} with Knight \cite{KnightPillaySteinhorn} shortly after
the original \cite{vdDries}.

\begin{theorem}[van den Dries-Miller]
\label{omin}
The structure $\rr _{\rm an,exp}$ defined this way is o-minimal. 
\end{theorem}

\subsection{Some definability}
The power of the $\rr _{\rm an,exp}$ structure comes from:
\begin{itemize}

\item The fact that it is o-minimal (Theorem \ref{omin});

\item Properties that can be deduced, using the o-minimality as well as more specific arguments; 

\item And the fact that many situations give rise to examples of $\rr _{\rm an,exp}$-definable subsets or functions. 

\end{itemize}

There is wide agreement that this structure 
(or more generally any similarly flexible
o-minimal structure) replies to Grothendieck's search for a {\bf tame topology}, that would be a theory of topology of sets and spaces that avoids the many pathological examples blocking the realization of our basic intuitions. 

Basic to the application to the André-Oort conjecture and the 
hyperbolic Ax-Lindemann theorem is the case of the Siegel upper half space. The defining conditions for  ${\mathbb H}_g$ are semialgebraic, so ${\mathbb H}_g$ is definable in $\rr _{\rm alg}$ and {\em a fortiori} in $\rr _{\rm an,exp}$. Peterzil-Starchenko showed the definability of the map to the quotient, in the following sense, with Klingler-Ullmo-Yafaev treating the case of Shimura varieties. 

\begin{theorem}
\label{pidef}
Let ${\mathscr F} \subset {\mathbb H}_g$ be a fundamental domain for the action of $SP(2g,\zz )$. Let $\overline {\mathscr F}$ be its closure. Then 
$\overline {\mathscr F}$ and $\overline {\mathscr F}$ are $\rr _{\rm an,exp}$-definable, and the projection map 
$$
\pi :  {\mathbb H}_g \rightarrow  {\mathcal A}_g
$$
when restricted to $\overline {\mathscr F}$ is an $\rr _{\rm an,exp}$-definable
function. 
\end{theorem}

One should beware that, because of its periodicity with respect to the action of 
$Sp(2g,\zz )$, the projection map $\pi$ itself can not be definable in an o-minimal structure. 

\subsection{The o-minimal Chow theorem}
The Peterzil-Starchenko theorem \cite{PeterzilStarchenko} provides a good
sample application illustrating the reach of definability in an o-minimal structure. 
 
A complex algebraic variety is a subset $X\subset \pp^n _{\cc}$ defined by algebraic equations. The projective space is covered by affine open subsets $\cc ^n \cong \rr^{2n}$ and the transition maps are algebraic, so there is a well-defined notion of definability on $\pp^n_{\cc}$ and $X$ is a definable subset (in any structure since all definability structures contain $\rr _{\rm alg}$). 

Assuming it is smooth, $X$ also has a structure of complex manifold. 
Serre's GAGA theorem tells us, among other things, that a closed subspace that is a holomorphic submanifold (that is to say, locally defined by holomorphic functions, or equivalently, having coordinate patches given by holomorphic functions) is in fact algebraic. This can fail badly for non-closed subsets. The first typical example is the graph of the exponential map $\{ (x,e^x)\}$ viewed inside $\cc^2$ and then in 
$\pp ^2_{\cc}$. We say that the exponential map $z\mapsto e^z$ on a complex variable has an {\em essential singularity} at $z=\infty$. 

One of the main characteristics of essential singularities is that they oscillate wildly in certain directions. For example, our complex exponential, when restricted to the imaginary axis $\sqrt{-1}\cdot \rr \subset \cc$, is given by sines and cosines:
$$
e^{\sqrt{-1} t} = \cos (t) + \sqrt{-1}\sin (t).
$$
Let us recall more generally that solutions of differential equations with {\em wild} (i.e. irregular) singularities always have ``Stokes directions'' in which such oscillatory behavior occurs. 

Peterzil and Starchenko \cite{PeterzilStarchenko} realized that this characteristic of essential singularities would allow to use definability in an o-minimal structure to rule out such a behavior. This leads to their {\em o-minimal Chow theorem:} 

\begin{theorem}[Peterzil-Starchenko]
\label{ominchow}
Suppose $X$ is a complex algebraic variety and $Y\subset X$ is a locally closed subset that is definable in an o-minimal structure. Suppose that locally $Y$ is a holomorphic submanifold or even a possibly singular subspace defined by holomorphic equations. If $Y$ is closed, or if the closure of $Y$ has no point at which it is a real manifold with nonempty real codimension $1$ boundary, then the Zariski closure $Y^{\rm Zar}$ is an algebraic subvariety containing $Y$ as a definable dense open subset. 
\end{theorem}

This has been generalized to the non-archimedean setting by Oswal \cite{Oswal}, whose paper also provides a brief review of many of our topics. 

\subsection{The o-minimal GAGA theorems} 

The Peterzil-Starchenko theorem is sometimes itself known as ``o-minimal GAGA'', but the classical GAGA (``Géométrie Algébrique et Géométrie Analytique'') theory goes farther. We recall that Serre formulated it mainly in terms of coherent sheaves on proper algebraic varieties, showing that analytic coherent sheaves are algebraic, together with similar statements for morphisms, section, and so forth. 

Bakker, Brunebarbe and Tsimerman \cite{BBT-Griffiths} develop various aspects generalizing Serre's GAGA theory to the definable setting. One of the major starting points is the definition of the category of coherent sheaves itself on a definable space. They prove that the functor from coherent sheaves on an algebraic variety to ones on its definable structure, is fully faithful, exact, with image closed under subquotients. Using the Peterzil-Starchenko theorem \ref{ominchow} this is done even for quasiprojective varieties, but on the other hand, the functor is not fully faithful. 

They go on to obtain an analogue of Theorem \ref{ominchow} saying that the image of a 
definable proper complex analytic morphism is algebraic. Here again we see that essential singularities are ruled out. 

They then apply this theory in their proof of the Griffiths conjecture, Theorem \ref{grifconj} below. 

\section{Ax-Lindemann theorems}
The hyperbolic Ax-Lindemann theorem \ref{axlindemann} is a generalization of several earlier theorems of the same type that were proven for definable subsets. It is instructive to look at how definability enters in here, and how some of the early proofs went. 

As we discussed in \S \ref{sec-limsum} above, the classical Lindemann-Weierstrass theorem was a first major generalization of the transcendence of $e$.  It gives algebraic independence of the exponentials $e^{a_1},\ldots , e^{a_n}$ over $\qq$. This was then extended by the Gelfond-Schneider theorem and  Baker's theorem giving algebraic independence of $\{ a_i^{b_i}\}$ (as well as quantitative bounds). These have a conjectural further generalization that remains an open problem: 

\begin{conjecture}[Schanuel conjecture]
For any complex numbers $z_1,\ldots , z_n$ linearly independent over $\qq$, the field 
$$
\qq ( z_1,\ldots , z_n,e^{z_1},\ldots , e^{z_n})
$$ 
has transcendence degree $\geq n$.
\end{conjecture}

An important milestone was James Ax's article \cite{Ax} introducing a functional version of the Schanuel conjecture and leading to a general collection of situations that are of main concern to us here, that may be roughly grouped under the terminology {\em functional transcendence}. The Ax-Schanuel theorem in \cite{Ax} says that if $f_1,\ldots , f_n\in t\cc[[t]]$ are $\qq$-linearly independent then the field $\cc (t,f_1,f_2,\ldots ,e^{f_1},e^{f_2},\ldots )$ has transcendence degree at least $n$ over $\cc ( t)$. 

This in turn implies the {\em functional Ax-Lindemann theorem} for the exponential function: given an irreducible algebraic subvariety $Y\subset \cc^n$ then the Zariski closure of its image by the exponential map $(z_1,\ldots , ) \mapsto 
(e^{2\pi i z_1}, \ldots  )$ is a translate of a subtorus.  Going in the other direction, if $Z\subset (\cc ^{\ast})^n$ is an algebraic subvariety then any irreducible algebraic subvariety of its inverse image by the exponential map 
is a translate of a  linear subspace. 

Klingler, Ullmo and Yafaev propose a very general formulation \cite{KlinglerUllmoYafaev}
of Ax-Lindemann properties in terms of subsets that are bi-algebraic for two distinct structures of complex algebraic variety on the same complex analytic space. Roughly speaking, the statement we expect is that the only bi-algebraic subsets are ones that we already know for geometric reasons.

We now get to Theorem \ref{axlindemann} and present a rough sketch of the core part of the proof \cite[proof of Theorem 6.1]{PilaTsimerman}. 

To keep with the notations of the article, we have an algebraic subvariety $V\subset {\mathcal A}_g$ and a real algebraic subset $W$ of the inverse image of $V$ in ${\mathbb H}_g$. After some reductions it is supposed that $W$ is maximal and $V$ is minimal, and furthermore $W$ is in fact a domain inside a complex algebraic subvariety of the Grassmanian that contains ${\mathbb H}_g$ as a bounded domain. 

Then $W$ has some algebraic curves $C\subset W$ and they show that $C$ must intersect the boundary $\partial {\mathbb H}_g$ in a sufficiently large subset (characterized by a volume estimate). 

The {\em monodromy group} $\Gamma \subset Sp(2g,\zz )$ means the monodromy group of the family of abelian varieties parametrized by $V$. It is infinite and acts on 
${\mathbb H}_g$, and this action extends to an action of the real group $G_0(\rr )$ 
(connected component of the real Zariski closure of $\Gamma$, basically the Mumford-Tate group of the VHS over $V$). 

Pila and Tsimerman define a subset $X\subset G_0(\rr )$ to be the set of group elements taking an appropriately chosen fundamental domain to a translate that touches $W$. We may think of this as being a subset of elements that approximately takes a basepoint in $W$ to another point in $W$, where the notion of ``approximately'' here is defined using the concept of how the fundamental domain touches the subspace. 

The study of volumes of curves $C\subset W$ near the boundary leads to the conclusion that $X$ has a relatively large number of points in $\Gamma$ of height bounded by $T$, namely growing as a strictly positive power of $T$. 

Then, $X$ is definable so the counting theorem \ref{counting} provides a decomposition of $X$ into algebraic pieces, one of which $W_1$ must account for a big enough number of points. An argument using maximality of $W$ yields a large number of points of 
$\Gamma$ leaving invariant $W$. Taking the Zariski closure leads, in turn, to the subgroup that will be used to define the Shimura subvariety needed to obtain a weakly special $S$. This concludes the proof sketch adapted to our present scope.

\section{Shimura varieties}
\label{sec-shimura}
The Siegel upper half-plane is one of the motivating examples that led to Griffiths' theory of period domains.
Next come the more general {\em Shimura varieties}. We have seen these in disguise, indeed the subspaces of ${\mathbb H}_g$ or ${\mathcal A}_g$ that we are calling ``special subvarieties'' defined by conditions on the Mumford-Tate group, are really Shimura subvarieties. 

The period domain of a Shimura variety takes the form of a quotient 
$$
{\mathscr H} = G(\rr ) / K(\rr )
$$
where $G$ is a real Lie group and $K\subset G$ is a maximal compact subgroup, such that $G$ is in a list of Lie groups of Hermitian type. This condition is designed to guarantee 
that ${\mathscr H}$ has a natural structure of complex manifold. 

One way of describing the Hermitian type condition is to say that the Lie algebra of $G$ should have a ``Hodge decomposition'' with three pieces
$$
{\mathfrak g} \otimes _{\rr} \cc  = {\mathfrak g}^{-1,1} \oplus 
{\mathfrak g}^{0,0}  \oplus {\mathfrak g}^{1,-1} 
$$
such that ${\mathfrak g}^{0,0} = {\mathfrak k}$ is the (complexified) Lie algebra of the maximal compact subgroup. The tangent space of ${\mathscr H} $ then gets a 
decomposition coming from the two pieces 
${\mathfrak g}^{-1,1}$ and ${\mathfrak g}^{1,-1}$, yielding the complex structure. 

The double quotient $G(\zz ) \backslash {\mathscr H}$ is a complex manifold (avoiding as we are doing, discussion of ``stacky points'' that come about from points in 
${\mathscr H}$ having a nontrivial finite stabilizer), that is the analogue for general Shimura varieties of the moduli space ${\mathcal A}_g$ of abelian varieties. 

The theory of the Satake-Baily-Borel compactification works also for 
$G(\zz ) \backslash {\mathscr H}$ and gives it a structure of quasiprojective variety. 

All of our discussion for the case of the Siegel upper half plane also works for the projection map 
$$
\pi : {\mathscr H} \rightarrow G(\zz ) \backslash {\mathscr H} .
$$
Namely, we have a hyperbolic Ax-Lindemann theorem and the André-Oort conjecture. These are due to Klingler, Ullmo, Yafaev, Mok, Pila, Shankar, Tsimerman  \cite{KlinglerUllmoYafaev,MokPilaTsimerman,PilaShankarTsimerman}. 

\begin{theorem}
\label{axshim}
The hyperbolic Ax-Lindemann theorem and the André-Oort conjecture hold for 
Shimura varieties. 
\end{theorem}

The Ax-Lindemann property was proven for Shimura varieties in \cite{KlinglerUllmoYafaev}. It is seen, in \cite{MokPilaTsimerman}, as a special case of a more general {\em Ax-Schanuel} property about algebraic relations between the coordinates on the Lie theory side and the algebraic geometry side. We'll state this more precisely below, as it extends to Griffiths period maps in greater generality. 

The proof of the André-Oort property (that we're calling MMMLAO) in \cite{PilaShankarTsimerman} is more than just a generalization of the case of ${\mathcal A}_g$, it invokes many new structures and techniques such as p-adic Hodge theory and the work of Esnault and Groechenig on rigid connections \cite{EsnaultGroechenig}. 

\section{More general period maps}

Let us look now more closely at the definition of the Griffiths period domains. This discussion is still not fully general: one could expand the setting to include general groups of Hodge type, but we'll stick to linear Hodge structures. Fresán's Bourbaki \cite{FresanBourbaki} provides an in-depth discussion and important references for what we discuss in this section. 

A {\em real Hodge structure of weight $k$} is a real vector space $V_{\rr}$
together with a polarization form $\langle \, , \, \rangle : V_{\rr} \times V_{\rr} \rightarrow \rr$ and a Hodge decomposition
$$
V_{\rr} \otimes _{\rr} \cc = \bigoplus _{p+q=k} V^{p,q}
$$
such that $\overline{V^{p,q}} = V^{q,p}$ and satisfying various compatibility conditions with the polarization form, notably a signed positive definiteness condition on the Hodge subspaces with the signs alternating in $(p,q)$. 

Given the ``Hodge numbers'' $h^{p,q}:= \dim (V^{p,q} )$, such real Hodge structures are ``all the same''. Technically speaking, the real algebraic group 
$$
G(\rr ) := Aut(V_{\rr}, \langle \, , \, \rangle )
$$
defined by the bilinear form (it is antisymmetric if $k$ odd, symmetric if $k$ even) acts transitively on the space of Hodge structures on the given space $V_{\rr}$. The isotropy group of elements that fix a given Hodge decomposition is compact, however it is not in general maximal compact. That is one of the really basic differences with the Siegel and Shimura cases. 

Thus, we can write down the {\em period domain}, the space of Hodge decompositions on a given $V_{\rr}, \langle \, , \, \rangle $, in the usual Lie theory picture as a quotient
$$
{\mathscr D} := G(\rr ) / K.
$$
The points of ${\mathscr D}$ parametrize exactly the Hodge decompositions compatible with $\langle \, , \, \rangle$. 

Part of the theory says that the decomposition is determined by the corresponding {\em Hodge filtration}, and the space of such filtrations is a complex algebraic flag variety known as the {\em compact dual} $\widecheck{\mathscr D}$. This complex dual is embedded as a closed algebraic subvariety of a product of Grassmanians. 

A $\zz$-Hodge structure  is just a real Hodge structure together with a lattice 
$V_{\zz} \subset V_{\rr}$ preserved (at least up to a rational constant)
by the polarization form. By ``lattice'' it means 
that $V_{\zz}\otimes _{\zz} \rr = V_{\rr}$. If the real Hodge structures are ``all the same'', that is not true for the $\zz$-Hodge structures which contain the numerical information of periods defining the position of the Hodge filtration. If $G(\zz )$ denotes
the group of transformations in $G(\rr )$ that preserve $V_{\zz }$ then the space parametrizing $\zz$-Hodge structure of the given type, up to isomorphism, is 
the quotient $G(\zz ) \backslash {\mathscr D}$. 

The theory looks so far the same as for Shimura varieties, generalizing the Siegel upper half space. There is a signiflcant difference. In the case of general Hodge structures, Griffiths defined the condition of {\em transversality} for holomorphic maps into 
$G(\zz ) \backslash {\mathscr D}$ or their liftings into 
${\mathscr D}$. The first point to note here is that there is indeed a complex structure, but the transversality condition also says that to first order the Hodge subspaces should vary only in the direction of their nearest neighbors. Thus, $V^{p,q}$ is allowed to tilt towards $V^{p-1,q+1}$ or $V^{p+1,q-1}$ in a first order deformation. 

It is easy to see that this condition is vacuous for Hodge structures of weight $k=1$, where there are only two subspaces. It is also vacuous in a few other cases, such as the case of weight $2$ with Hodge numbers 
$$
h^{0,2} = 1, h^{1,1} , h^{2,0} = 1.
$$
That occurs for example when we look at K3 surfaces where $h^{1,1} = 20$. The cases where the Griffiths transversality condition holds automatically correspond to the cases where $G(\rr )$ is of Hermitian type, and hence when the period domain is a Shimura variety. 

In general, though, the transversality condition is an important condition. The name comes from the fact that it isolates tangent directions that are transverse to the compact subvarieties of ${\mathscr D}$ obtained by looking at $U(\rr ) / K(\rr )$ where 
$U(\rr )\subset G(\rr )$ is a maximal compact subgroup. 

As we saw in Section \ref{sec-traag}, given a family of varieties $f:X\rightarrow S$, the family of Hodge structures of the cohomology groups $H^k(X_s)$ becomes a family of $\zz$-Hodge structures, hence it determines a {\em period map}
$$
\Phi :S \rightarrow G(\zz ) \backslash {\mathscr D}
\;\; \mbox{ lifting to } \;\; 
\widetilde{\Phi} : \widetilde{S} \rightarrow {\mathscr D}
$$
with its lifting to the universal cover. The map $\widetilde{\Phi}$ is holomorphic and satisfies the Griffiths transversality condition. 

The transversality condition defines a collection of subspaces in the tangent bundle of ${\mathscr D}$, that is sometimes called a {\em distribution} (not to be confused with distributions as generalized functions in analysis). It alows to recover some of the classical properties of Hermitian symmetric domains such as the Siegel upper half plane. For example, there is a natural metric having negative curvature in the directions of the horizontal distribution. This partial negativity of the curvature is sufficient to obtain many  local analytic properties generalizing those of the Hermitian symmetric case. 

The horizontal 
distribution is however {\bf not integrable}, that is to say it does not come from a foliation. It is more closely analogous to a ``contact structure'', and in the general case one can move from any point of ${\mathscr D}$ to any other point along zig-zag paths traveling tangent to the horizontal distribution. For this reason, we can not think of 
$G(\zz ) \backslash {\mathscr D}$ as being completely analogous to the moduli space of abelian varieties ${\mathcal A}_g$.

\subsection{Definability and the Ax-Schanuel theorem}

The works of a group of people including 
Bakker, Brunebarbe, Klingler, Mok, Pila, Shankar, Tsimerman, Ullmo and Yafaev led to the formulation and proof for general period domains of several of the main theorems that have been seen above for the case of ${\mathcal A}_g$.

A first main point is the definability of Griffiths period maps on fundamental domains \cite{BKT}. It is interesting to note that the potential for application of definability to local singularities such as appear in Hodge theory, was forseen early on with the paper of Schmid and Vilonen \cite{SchmidVilonen} that is referred to in Wilkies' Bourbaki notes \cite{WilkieBourbaki}. 

For the hyperbolic Ax-Lindemann theorem, Klingler formulated a more general statement for period domains, the {\em Ax-Schanuel property} \cite{Klingler}. To motivate this, we notice first that, given a definable map, its inverse map is also definable, just by taking the graph and 
flipping the two coordinate axes. 

In particular, in the Siegel modular case, the projection map 
$\pi : {\mathbb H}_g \rightarrow {\mathcal A}_g$ may be considered in the other direction, yielding exactly the period map for the universal family of abelian 
varieties\footnote{It may be worth repeating once again that we are ignoring stacky aspects at points where the group action has an isotropy group, corresponding to abelian varieties with automorphisms.}
$$
\widetilde{\Phi} : {\mathcal A}_g \rightarrow Sp(2g,\zz ) \backslash 
{\mathbb H}_g .
$$
The Ax-Schanuel property aims to be symmetric with respect to the two sides, the Lie-theory side and the algebraic geometry side. 

Suppose we are given a variety $X$ with a map $\Phi : X \rightarrow 
{\mathscr D}/\Gamma $ (here denoting by $\Gamma$ a general discrete group of transformations and writing the quotient now on the right). Let $\widetilde{X}$ denote the universal covering. We obtain the subset 
$$
\widetilde{Corr}_{\Phi} := X\times _{{\mathscr D}/\Gamma} {\mathscr D}
\subset X \times _{{\mathscr D}/\Gamma} {\mathscr D}
$$
which is basically the lifted graph of $\Phi$. More precisely it  is the correspondence
consisting of the pairs
of points $(x,y)$ satisfying the relationship that $y\in {\mathscr D}$ is a representative for $\Phi (x)$.

The {\bf Ax-Schanuel property} was formulated by Klingler \cite{Klingler} for general VHS. It states that this correspondence does not contain any subspaces satisfying any unexpected algebraic equations, with respect to the algebraic structures of $X$ as an algebraic variety, and ${\mathscr D}$ as an open subset of the 
compact dual flag variety $\widecheck{\mathscr D}$. 
By ``unexpected'', it means equations not implied by relationships between Mumford-Tate groups as would be consequences of the Hodge conjecture. 
We refer to \cite{Klingler,MokPilaTsimerman,BT} for the technical statement involving codimension bounds. 

The Ax-Schanuel property was proven for Shimura varieties by  Mok, Pila and Tsimerman \cite{MokPilaTsimerman} and for general pure VHS by Bakker and Tsimerman \cite{BT}. A version exists for VMHS that goes beyond our scope here. 
The Ax-Lindemann property is a consequence. 

\begin{theorem}[Bakker-Tsimerman \cite{BT}: Ax-Schanuel for VHS]
If $\Phi$ is the period map of a $\zz$-variation of Hodge structures over a
quasiprojective variety, then:
\begin{itemize}

\item
The correspondence $\widetilde{Corr}_{\Phi}$ is definable in the 
$\rr _{\rm an,exp}$ o-minimal structure; and

\item
This correspondence satisfies the Ax-Schanuel property. 
\end{itemize}
\end{theorem}

This proof cointains also the Ax-Schanuel result for Shimura varieties
\cite{MokPilaTsimerman} that we stated previously, leading to Theorem \ref{axshim}. Lawrence and Venkatesh show that this result in the complex setting leads to a $p$-adic analogue \cite{LawrenceVenkatesh}.

\subsection{The Griffiths period image conjecture}

The Ax-Schanuel theorem shows that the theory of period maps can enter into a theory similar to what was developed for ${\mathcal A}_g$. However, for general period maps the link between the Lie theory side and the algebraic geometry side remains less strong. It has indeed been a longstanding motivating question in the theory of variations of Hodge structure, to provide as much algebraic geometry as possible to the period domain. 

Because of the horizontality condition, non-integrability of the horizontal distribution, and the compact vertical subvarieties, there is always going to be a certain disparity. Griffiths formulated his conjecture on the images of period maps in order to reduce the gap. 

Much of the study of local properties of variations of Hodge structure near singular divisors, by Griffiths, Schmid, Cattani, Kaplan, Kashiwara, Kawai was motivated, in part, by being able to gain some control over the boundary behavior of period maps. 

It turns out that these studies did mesh perfectly well with the viewpoint of o-minimal geometry, the ``tame topology'' perspective needed to put the required information together. The Griffiths conjecture is now a theorem by Tsimerman with Bakker and Brunebarbe \cite{BBT-Griffiths}: 

\begin{theorem}[Bakker-Brunebarbe-Tsimerman: Griffiths conjecture]
\label{grifconj}
If $X$ is a quasiprojective variety and $\Phi : X \rightarrow {\mathscr D} / \Gamma$
is the period map of a $\zz$-variation of Hodge structure, then the image $\Phi (X)$ has the structure of a quasiprojective variety such that $\Phi$ becomes an algebraic map. 
\end{theorem}

This theorem opens up the way towards a future study of the arithmetic and algebraic geometry of these period images. We mention here Griffiths' original motivation: to use the period images in order to construct moduli spaces for varieties, whenever the (possibly infinitesimal) Torelli property is known. 

Indeed, the Torelli property is the statement that the moduli point of an algebraic variety is accurately reflected in its corresponding points in the Griffiths period domains, i.e. by the Hodge structures on its cohomology group. The infinitesimal version is the same statement but just for deformations of a given variety. Torelli theorems abound, for curves, hypersurfaces, and many other types of varieties. When this property does hold, then the image of the period map will serve as a construction of the moduli space for the varieties in the given deformation class. 

\subsection{The Shafarevich conjecture for quasi-projective varieties} 
Applications of Tsimerman's viewpoint to questions in algebraic geometry keep on arriving. We are not able to delve into the details of the most advanced among these, but let us look at one: the Shafarevich conjecture for linear representations of fundamental groups of quasiprojective varieties. 

The general Shafarevich conjecture predicts that the universal covering of an algebraic variety $X$ (over $\cc$) is {\em holomorphically convex}, saying that the holomorphic functions exhaust it in the sense that an unbounded sequence of points has a holomorphic function that is unbounded on it. Associated to this concept is the {\em Shafarevich morphism} that contracts subvarieties of $X$ whose fundamental groups have finite images in $\pi _1(X)$. Although the Shafarevich morphism has been constructed, following ideas of Kollár and others, we do not in general know anything about the universal covering of its base space. 

This changes when we transform the question using a representation (or collection of representations) $\rho : \pi _1(X) \rightarrow G$. We may then take the covering of $X$ generated by the representations (i.e. corresponding to the subgroup of $\pi _1(X)$ that is the intersection of their kernels) and ask for it to be holomorphically convex. The {\em relative Shafarevich morphism} is a map contracting the subvarieties of $X$ whose fundamental groups have finite images under $\rho$. 

In some cases, such as for our main example ${\mathcal A}_g$, there is a faithful representation of the fundamental group. We may note that in this example the universal covering ${\mathbb H}_g$ is a bounded domain so it is easily seen to be holomorphically convex (fractions involving the linear functions of the matrix coefficients of the period matrix coordinate for 
${\mathbb H}_g$ will do the trick). On the other hand, it is known after examples due originally to Toledo, that there are smooth projective varieties whose fundamental group is not residually finite, and therefore can not have a faithful collection of linear representations. Thus, the relative Shafarevich conjecture is a middle ground. 

Work of many mathematicians, culminating in the paper of Eyssidieux, Katzarkov, Pantev and Ramachandran \cite{EKPR}, led to the proof of the relative Shafarevich conjecture for smooth projective varieties. 

The quasiprojective case required new ideas, and Tsimerman with co-authors Bakker and Brunebarbe \cite{BBT-Sha} were able to supply these, including notably the concept of definability in an o-minimal structure. They used ingredients such as Mochizuki's theory \cite{Mochizuki} and made further developments in the theory of the Deligne-Hitchin twistor space for quasiprojective varieties. See \cite{Bakker} for a survey. 

\begin{theorem}[Quasi-projective Shafarevich conjecture: Bakker-Brunebarbe-Tsimerman] 
Suppose $X$ is a quasiprojective variety over $\cc$, and suppose given a collection of representations $\rho : \pi _1(X)\rightarrow G$ to linear algebraic groups. Suppose that this collection satisfies a condition of defining a nice subset of the moduli space of representations (satisfied in particular if we take all representations of a given rank). Then the Shafarevich conjecture holds relative to $\{ \rho \}$, that is to say the covering of $X$ defined by the intersection of kernels of the $\rho$'s is holomorphically convex. 
\end{theorem}

This remarkable journey has taken us from the o-minimal definability theory of special points and special subvarieties of modular varieties naturally occuring in arithmetic geometry, towards the most general structures coming out of moduli theory in the algebraic geometry of higher dimensional varieties. Jacob Tsimerman's works have taken us on a wild ride through tame topology.

\section*{Acknowledgments.}
I would like to thank Jacob Tsimerman and Jonathan Pila for their help and encouragement in the preparation of this article, and to thank all those who have contributed to my understanding of this subject.

No dedicated AI tools were used for this article. However, the AI-mode responses to Google search queries were often useful. 

This work was supported by the 
European Research Council through the Horizon ERC Synergy grant 
101167526 (MALINCA) and by the ANR program IsoMoDyn 
(ANR-25-CE40-1360).

\end{document}